\documentclass[12pt]{article}

\usepackage{amssymb, amsmath, amsfonts, amsthm, mathrsfs, latexsym, epsfig,  psfrag, graphicx}
\usepackage[all]{xy}
\xyoption{v2} \xyoption{knot}\xyoption{arc} \xyoption{poly}
\xyoption{curve}

\newcommand{\bbc}{\mathbb{C}}

\newcommand{\aut}{\mathrm{Aut}}

\usepackage{amsmath,amsthm,amsfonts,amssymb,amscd}
\def\qed{{\unskip\nobreak\hfil\penalty50
\hskip2em\hbox{}\nobreak\hfil$\square$
\parfillskip=0pt \finalhyphendemerits=0\par}\medskip}

\def\prf{\trivlist \item[\hskip \labelsep{\bf Proof\ }]}

\def\Ad{{\mathrm {Ad}}}

\def\C{{\Bbb {C}}}
\def\End{{\mathrm {End}}}

\def\Hom{{\mathrm {Hom}}}

\def\dim{{\mathrm {dim}}}

\def\Ad{{\mathrm {Ad}}}

\def\End{{\mathrm {End}}}
\def\Hom{{\mathrm {Hom}}}

\def\dim{{\mathrm {dim}}}

\newtheorem{theorem}{Theorem}[section]
\newtheorem{lemma}[theorem]{Lemma}
\newtheorem{conjecture}[theorem]{Conjecture}
\newtheorem{corollary}[theorem]{Corollary}
\newtheorem{definition}[theorem]{Definition}

\newtheorem{proposition}[theorem]{Proposition}
\newtheorem{remark}[theorem]{Remark}

\def\Hom{{\mathrm{Hom}}}

\def\A{{\cal A}}
\def\L{{\cal L}}
\def\F{{\cal F}}

\def\B{{\cal B}}

\def\Z{{\mathbb Z}}

\def\C{{\mathbb C}}

\renewcommand{\qed}{\ \hfill $\blacksquare$}

\newcommand{\bdef}{\begin{definition}}
\newcommand{\blem}{\begin{lemma}}
\newcommand{\bprop}{\begin{proposition}}
\newcommand{\bthm}{\begin{theorem}}
\newcommand{\bcor}{\begin{corollary}}
\newcommand{\bconj}{\begin{conjecture}}
\newcommand{\ede}{\end{definition}}
\newcommand{\elem}{\end{lemma}}
\newcommand{\eprop}{\end{proposition}}
\newcommand{\ethm}{\end{theorem}}
\newcommand{\ecor}{\end{corollary}}
\newcommand{\econj}{\end{conjecture}}
\newcommand{\brem}{\begin{remark}}
\newcommand{\erem}{\end{remark}}

\newcommand{\ba}{\begin{array}}
\newcommand{\ea}{\end{array}}
\newcommand{\bea}{\begin{eqnarray}}

\title{\huge On a subfactor generalization of Wall's conjecture \\}
\author{
{\sc  Robert Guralnick}\footnote{Supported in part by NSF grants DMS- 0653873 and 1001962.}\\
Department of Mathematics\\
University of Southern California\\
Los Angeles, California 90089-2532\\
E-mail: {\tt guralnic@usc.edu}\\
{}\\
{\sc Feng Xu}\footnote{Supported in part by NSF grant DMS-0800521 and an academic senate grant from UCR.}\\
Department of Mathematics\\
University of California at Riverside\\
Riverside, CA 92521\\
E-mail: {\tt xufeng@math.ucr.edu}}

\begin{document}
\date{}
\maketitle

\begin{abstract}
In this paper we discuss a conjecture on intermediate subfactors
which is a generalization of Wall's conjecture from the theory of
finite groups.  We explore special cases of this conjecture and
present supporting evidence.  In particular we prove special cases
of this conjecture related to some finite dimensional Kac Algebras
of Izumi-Kosaki type which include relative version of Wall's
conjecture for solvable groups.

\end{abstract}

\newpage

\section{Introduction}

Let $M$ be a factor represented on a Hilbert space and $N$ a
subfactor of $M$ which is irreducible, i.e.,$N'\cap M= \C$. Let $K$
be an intermediate von Neumann subalgebra for the inclusion
$N\subset M.$ Note that $K'\cap K\subset N'\cap M = \C,$ $K$ is
automatically a factor. Hence the set of all intermediate subfactors
for $N\subset M$ forms a lattice under two natural operations
$\wedge$ and $\vee$ defined by:
\[
K_1\wedge K_2= K_1\cap K_2, K_1\vee K_2= (K_1\cup K_2)''.
\]
The commutant map $K\rightarrow K'$ maps an intermediate subfactor
$N\subset K\subset M$ to $M'\subset K'\subset N'.$ This map
exchanges the two natural operations defined above.

Let $M\subset M_1$ be the Jones basic construction of $N\subset M.$
Then  $M\subset M_1$ is canonically isomorphic to $M'\subset N'$,
and the lattice of intermediate subfactors for $N\subset M$  is
related to the lattice of intermediate subfactors for $M\subset M_1$
by the commutant map defined as above.

 Let $G_1$ be a group and $G_2$ be a subgroup of $G_1$. An
interval sublattice $[G_1/G_2]$ is the lattice formed by all
intermediate subgroups $K, G_2\subseteq K\subseteq G_1.$

By cross product construction and Galois correspondence,  every
interval sublattice of finite groups can be realized as intermediate
subfactor lattice of finite index. Hence the study of intermediate
subfactor lattice of finite index is a natural generalization of the
study of interval sublattice of finite groups. The study of
intermediate subfactors has been very active in recent years(cf.
\cite{BJ},\cite{GJ}, \cite{JXu},\cite{Jh2}, \cite{ILP},
\cite{Longo4}, \cite{Wat} and \cite{ SW} for only a partial list).
By a result of S. Popa (cf. \cite{Popa}), if a subfactor $N\subset
M$ is irreducible and has finite index, then the set of intermediate
subfactors between $N$ and  $M$ is finite. This result was also
independently proved by Y. Watatani (cf. \cite{Wat}). In \cite{Wat},
Y. Watatani investigated the question of which finite lattices can
be realized as intermediate subfactor lattices. Related questions
were further studied by P. Grossman and V. F. R. Jones in \cite{GJ}
under certain conditions. As emphasized in \cite{GJ}, even for a
lattice consisting of six elements with shape a hexagon, it
is not clear if it can be realized as intermediate subfactor lattice
with finite index. This question has been solved recently by M.
Aschbacher in \cite{Asch} among other things. In \cite{Asch}, M.
Aschbacher  constructed a finite group $G_1$ with a subgroup $G_2$ such
that the interval sublattice $[G_1/G_2]$ is a hexagon. The lattices
that appear in \cite{GJ,Wat,Asch} can all be realized as interval
sublattice of finite groups. There are a number of old problems
about interval sublattice of finite groups. It is therefore a
natural programme to investigate if these old problems have any
generalizations to subfactor setting. The hope is that maybe
subfactor theory can provide new perspective on these old
problems.\par In \cite{Xl} we consider the problem whether the very
simple lattice $M_n$ consisting of a largest, a smallest and $n$
pairwise incomparable elements can be realized as subfactor lattice.
We showed in \cite{Xl} all $M_{2n}$ are realized as the lattice of
intermediate subfactors of a pair of hyperfinite type $III_1$
factors with finite depth. Since it is conjectured that infinitely
many $M_{2n}$ can not be realized as interval sublattices of finite
groups (cf. \cite{Balu} and \cite{Palfy}), our result shows that if
one is looking for obstructions for realizing finite lattice as
lattice of intermediate subfactors with finite index, then the
obstruction is very different from what one may find in finite group
theory.\par

In 1961 G. E. Wall conjectured that the number of maximal subgroups
of a finite group $G$ is less than $|G|$, the order of $G$ (cf.
\cite{wall61}). In the same paper he proved his conjecture when $G$
is solvable. See \cite{lie} for more recent result on Wall's
conjecture.

Wall's conjecture can be naturally generalized to a conjecture about
maximal elements in the lattice of intermediate subfactors. What we
mean by maximal elements are those subfactors $K\neq M,N $ with the
property that if $K_1$ is an intermediate subfactor and $K\subset
K_1,$ then $K_1=M$ or $K.$ Minimal elements are defined similarly
where $N$ is not considered as an minimal element. When $M$ is the
cross product of $N$ by a finite group $G$, the maximal elements
correspond to maximal subgroups of $G,$ and the order of $G$ is the
dimension of second higher relative commutant. Hence a natural
generalization of Wall's conjecture as proposed in \cite{X1} is the
following:
\begin{conjecture}\label{wall}
Let $N\subset M$ be an irreducible subfactor with finite index. Then
the number of maximal  intermediate subfactors is less than
dimension of $N'\cap M_1$ (the dimension of second higher relative
commutant of $N\subset M$).
\end{conjecture}
We note that since maximal intermediate subfactors in $N\subset M$
correspond to minimal intermediate subfactors in $M\subset M_1,$ and
the dimension of second higher relative commutant remains the same,
the conjecture  is equivalent to a similar conjecture as above with
maximal replaced by minimal.\par In \cite{X1}, Conjecture \ref{wall}
is verified for subfactors coming from certain conformal field
theories. These are subfactors not related to groups in general. In
this paper we consider those subfactors which are more closely
related to groups and more generally Hopf algebras.\par

If we take $N$ and $M$ to be cross products of a factor $P$ by $H$
and $G$ with $H$ a subgroup of $G$, then  the minimal version of
conjecture \ref{wall} in this case states that the number of minimal
subgroups of $G$ which strictly contain $H$ is less than the number
of double cosets of $H$ in $G$. This follows from simple counting
argument. The nontrivial case is the maximal version of the above
conjecture. In this case it gives a generalization of Wall's
conjecture which we call relative version of Wall's conjecture. The
relative version of Wall's conjecture states that the number of
maximal subgroups of $G$ strictly containing a subgroup $H$ is less
than the number of double cosets of $H$ in $G.$  As a simple example
when this can be proved, consider $G=H\times H,$ and $D \cong H$ is a
diagonal subgroup of $G.$ Then the set of maximal subgroups of $G$
containing $D$ are in one to one correspondence with the set of
maximal normal subgroups of $H,$ and it is easy to check that the set
of maximal normal subgroups has cardinality less than the number of
irreducible representations of $H.$ On the other hand  the number of
double cosets of $H$ in $G$ is the same as the number of conjugacy
classes of $H,$ and this is the same as the number of irreducible
representations of $H.$ So we have proved the relative version of
Wall's conjecture in this case.

In \S\ref{relativesection} we will prove this relative version of
Wall's conjecture for  $G$ solvable.   We
will present two proofs. The first  proof is motivated by an idea of
V. F. R. Jones which is to seek linear independent vectors
associated with minimal subfactors in the space of second higher
relative commutant. This proof is indirect but we hope that the idea
will prove to be useful for more general case. We formulate a
conjecture for general subfactors (cf. Conjecture \ref{mod}) which
is stronger than Conjecture \ref{wall}, and for solvable groups this
conjecture is proved in \cite{X1}. Here we modify the proof in
\cite{X1} to prove a linear independence result (cf. Th. \ref{sol}),
and this result implies the relative Wall conjecture for solvable groups.
 The second proof is a
more direct proof using properties of maximal subgroups of solvable groups.
\par The cross
product by finite group subfactor  is a special case of depth $2$
subfactor. If we take $N\subset M$ to be depth 2, by
\cite{David},\cite{NV} such a subfactor comes from cross product by
a finite dimensional *-Hopf algebra or Kac algebra $\A$. By
\cite{ILP} or \cite{NV} the intermediate subfactors are in one to
one correspondence to the set of left (or right) coideals of $\A$.
Then Conjecture \ref{wall} states that the number of maximal (resp.
minimal) right coideals of $\A$ is less than the dimension of $\A$.
In \S\ref{iksection} we will prove this for the case of Kac algebras
$\A$ of Izumi-Kosaki type with solvable groups as considered in
\cite{IK}. We also prove Conjecture \ref{wall} for the intermediate
subfactors of Izumi-Kosaki type with solvable groups as considered
in \cite{IK} which are not necessarily of depth 2 (cf. Th.
\ref{relative2}). Th. \ref{relative2} generalizes Th.
\ref{relativeth}. It is interesting to note that the same type of
first cohomology problem encountered in Remark \ref{coh} also
appears here but in a different way and solvability is once again
used to ensure that the first cohomology group is trivial (cf. Lemma
\ref{count1} and Lemma \ref{count2}).\par We note that recently
lattices of intermediate for other types of Kac algebras have been
obtained in \cite{DT}. Our conjecture can be verified in the
examples of \cite{DT} where complete lattice of intermediate
subfactors are determined. The maximal (or minimal) coideals are
very few compared with the dimension of the Kac algebra in these
examples of \cite{DT}.
\par

In \S\ref{tensorsection} we first present a  lemma which bounds the
number of maximal subgroups of a group $X\times Y$ which does not
contain either $X$ nor $Y.$ This lemma gives a proof of Wall's
conjecture for $X\times Y$ assuming that Wall's conjecture is true
for $X$ and $Y.$ We then propose a natural conjecture about tensor
products of subfactors.

At the end of this introduction let us consider a fusion algebra
version of Conjecture \ref{wall}. Let $\rho_i\in \End(M), i=1,...,n$
be a finite system of irreducible sectors of a properly infinite
factor $M$ which is closed under fusion. Consider the Longo-Rehren
subfactor associated with such a system (cf. \cite{LR}). By
\cite{I}, the intermediate subfactors are in one to one
correspondence with the fusion subalgebras which are generated by a
subset of simple objects $\rho_i$ , and Conjecture \ref{wall} states
that the number of such maximal fusion subalgebras is bounded by $n$
which is the number of simple objects. This motivates us to make the
following conjecture:
\begin{conjecture}\label{fusioncase}
Let $\F$ be a finite dimensional semisimple fusion algebra with $n$
simple objects. Then  the number of maximal fusion subalgebras which
are generated by a subset of the simple objects of $\F$ is less than
$n$.
\end{conjecture}
If we take $\F$ to be the group algebra of $G,$ then Conjecture
\ref{fusioncase} is equivalent to Wall's conjecture.\par

If we take $\F$ to be the fusion algebra of representations of a
finite group $G$, then the maximal fusion subalgebras are in one to
one correspondence to minimal normal subgroups of $G,$ and the
number of such subgroups are less than the number of conjugacy
classes of $G$, which is the same as the number of simple objects of
$\F$. This is a special case of a more general result of   D.
Nikshych and V. Ostrik, who prove that Conjecture \ref{fusioncase}
is true for commutative $\F$ \cite{NO}.

The second named author (F.X.) would like to thank  Prof. V. F. R.
Jones for his encouragement and useful comments on conjecture
\ref{wall} which inspired our first proof presented in
\S\ref{relativesection}, and for communications on numerical
evidence supporting the relative version of Wall's conjecture. F.X.
would also like to thank Prof. D. Bisch for invitation to a
conference on subfactors and fusion categories in Nashville where
some of the results of this paper were discussed, and Professors
Marie-Claude David,  D. Nikshych and V. Ostrik for useful
communications.
\section{Relative version of Wall's
conjecture for solvable groups}\label{relativesection}

In this section we will prove Theorem \ref{relativeth}, which
confirms the relative version of Wall's conjecture for solvable
groups. We will give two proofs of this result. The first proof is
motivated by the following  conjecture, formulated as Conjecture A.1
in \cite{X1}, which can be stated for general subfactors:

\begin{conjecture}\label{mod}
Let $N\subset M$ be an irreducible subfactor with finite Jones
index, and let $P_i, 1\leq i\leq n$ be the set of minimal
intermediate subfactors. Denote by $e_i\in N'\cap M_1, 1\leq i\leq
n$ the Jones projections $e_i$ from $M$ onto $P_i$ and $e_N$ the
Jones projections $e_N$ from $M$ onto $N.$ Then there are vectors
$\xi_i, \xi\in N'\cap M_1$ such that $e_i \xi_i=\xi_i,1\leq i\leq n,
e_N\xi=\xi,$ and $\xi_i, 1\leq i\leq n,\xi$ are linearly
independent.
\end{conjecture}

\begin{remark}
We note that unlike conjecture \ref{wall}, the conjecture above
makes use of the algebra structure of $N'\cap M_1$ and therefore
does not immediately imply the dual version or if one replaces
minimal by maximal.
\end{remark}

By definition conjecture \ref{mod} implies conjecture \ref{wall}. In
the case of subfactors from groups,  it is easy to check that
conjecture \ref{mod} is equivalent to:

\begin{conjecture}\label{mod2'}
Let $K_i, 1\leq i\leq n$ be a set of maximal subgroups of $G.$ Then
there are vectors $\xi_i\in l(G), 1\leq i\leq n$ such that $
e_G\xi_i=0,$ $\xi_i$ are $K_i$ invariant and linearly independent.
\end{conjecture}

This conjecture is proved in \cite{X1}  when $G$ is solvable. It
turns out a modification of the proof presented in \cite{X1} gives a
proof of a stronger statement. Let us make the following stronger
conjecture.   First we need to introduce some notation.
If $H$ is a subgroup of $G$, let $\ell(H)=\ell(G,H)$ be the permutation
module $\bbc_H^G$.  Let $\ell_0(H)$ denote the hyperplane of
weight zero vectors in $\ell(H)$ (i.e. the complement to the $1$-dimensional
$G$-fixed space on $\ell(H)$).

\begin{conjecture}\label{mod2}
Let $K_i, 1\leq i\leq n$ be a set of maximal subgroups of $G.$
Set $H = \cap K_i$.  Then
there are vectors $\xi_i\in \ell_0(H), 1\leq i\leq n$  that
 are $K_i$-invariant and linearly independent.   In particular,
 this implies that $n \le \dim \ell_0(H)^H < |H / G \backslash H|$.
\end{conjecture}

We will prove Conjecture \ref{mod2} for solvable groups by modifying
the arguments of \cite{X1}.       We begin with some preparations
that hold for all finite groups.

\begin{lemma}\label{prim}
Suppose that $K_1, \ldots, K_n$ are conjugate maximal subgroups
of the finite group $G$.   Then Conjecture \ref{mod2} holds for
$\{K_1, \ldots, K_n\}$.
 \end{lemma}

\prf  Set $K=K_1$.  If $n=1$, the result is obvious (in particular, if $K$ is normal in $G$).  So assume that $n > 1$.  Let $H = \cap K_i$.
Of course, $\ell_0(K)$  is a submodule of $\ell_0(H)$.    Let $K_1, \ldots, K_m, m \ge n$ be the set of
all conjugates of $K$.  Since $K$ is not normal in $G$,  $K$ is self normalizing whence
if we choose a permutation basis $\{ v_i   | 1 \le i \le m\}$ for $\ell(K)$, then the stabilizers of the $v_i$ are precisely the $K_i$.
If $m > n$, then the vectors $v_i - v_0 \in \ell_0(K), 1 \le i \le n$ are clearly linearly independent (here $v_0 = \sum v_i$ is fixed
by $G$).   So it suffices to assume that $m=n$ and so in particular, $H$ is normal in $G$.  So we may assume
that $H=1$.   Let $V$ be a nontrivial irreducible submodule of $\ell_0(K)$.   Then $K$ does not act trivially on $V$.
Note that, by Frobenius reciprocity,  the multiplicity of $V$ in $\ell_0(K)$ is precisely $\dim V^K  < \dim V$.  Of course
$\dim V$ is the multiplicity of $V$ in $\ell_0(H)$.
  Thus,  $\ell_0(K) \oplus V$ is a submodule of $\ell_0(H)$.
  Now choose vectors $v_i - v_0, 1 \le i < m$  as above and $w_m$ any fixed vector of $K_m$ in $V$.  These are
obviously linearly independent.
 \qed

 \par
 Next we prove a reduction theorem for Conjecture \ref{mod2}.   Note that the reduction depends on
 the existence of the vectors and not just on cardinality.

\begin{lemma}\label{induction}   Let $\mathcal{S}$ be a family of finite simple groups.
Let $\mathcal{F(S)}$ denote the family of all finite groups with all composition factors
in $\mathcal{S}$.   Let $K_1, \ldots, K_n$ be maximal subgroups of the finite group
$G$ in $\mathcal{F(S)}$ and assume that Conjecture \ref{mod2} fails with $n|G|$ minimal.
Then each $K_i$ has trivial core in $G$.  In particular, $G$ is a primitive permutation group.
\end{lemma}

\prf   Suppose that $N$ is a nontrivial normal subgroup of $G$ contained in $K_1$.
Set $H = \cap K_i$.
If each $K_i$ contains $N$,  then $\ell_0(H)$ is a $G/N$-module and so
$G/N, K_1/N, \ldots, K_n/N$ would give a counterexample to the conjecture.

Reorder the $K_i$ so that $N \le K_i$  if and only if $i \le s < n$.   Note that $NK_j=G$ for
$j > s$.   By the minimality of $|G|n$,  we can choose  $v_1, \ldots, v_n  \in \ell_0(H)$
with $K_jv_j=v_j$ for  all $j$ such that $\{v_1, \ldots, v_s\}$ and $\{v_{s+1}, \ldots, v_n\}$
are linearly independent.    It thus suffices to show that spans of $v_1, \ldots, v_s$ and
$v_{s+1}, \ldots, v_n$ have trivial intersection.  Suppose that $u$ is in this intersection.

Since $e_Ne_{K_j}=e_G$ for $j > s$ (since $G=NK_j$), it follows that $0 = e_Gv_j = e_Ne_{K_j}v_j=e_Nv_j$
for $j > s$.  Thus, $e_Nu=0$.  Since $N$ fixes $v_i, i \le s$, it follows that $e_Nu=u$.  Thus, $u=0$
and the result follows.
\qed

\begin{theorem}\label{sol}
Conjecture \ref{mod2} is true for $G$ solvable.
\end{theorem}

\prf   Consider a counterexample with $|G| n$ minimal.  By Lemma
\ref{induction}, none of the $K_i$ contain a normal subgroup.  It
follows that $G$ is a solvable primitive permutation group, whence
$G = AK$ where $A$ is elementary abelian and $K$ acts irreducibly on
$A$. In particular, any maximal subgroup of $G$ either contains $A$
or is a complement to $A$. Since the core of each $K_i$ is trivial,
 $G=AK_i$ for each $i$.
Since $G$ is solvable, $H^1(K,A)=0$, whence all of the $K_i$ are conjugate.
Now apply Lemma \ref{prim} to complete the proof.
\qed

\begin{remark}\label{coh}  As we have seen, a minimal counterexample
to Conjecture 2.4 would be a
primitive permutation groups, and the set of maximal subgroups must all have
trivial core. Such groups are classified by Aschbacher-O'Nan-Scott
theorem (cf. \S 4 of \cite{NS}). The first case is when  $G$ is the
semidirect product of an elementary abelian group $V$ by $K_1,$ and
the action of $K_1$ on $V$ is irreducible. When $G$ is not solvable,
maximal subgroups $K$ of $G$ with trivial core are not conjugates of
$K_1$, and our proof as above does not work. Such maximal subgroups
are related to the first cohomology of $K_1$ with coefficients in
$V,$ and conjecture \ref{mod2} implies that the order of this
cohomology is less than $|K_1|$ (cf. Question 12.2 of \cite{lub}).
Unfortunately even though it is believed that the order of this
cohomology is small (cf. \cite{Gu}), the bound $|K_1|$ has not been
achieved yet.
\end{remark}

We give a second proof of Conjecture \ref{mod2} for solvable groups which
is not inductive.

Let $G$ be a solvable group.  Let  $H \le G$ and let   $K_1, \ldots, K_r$
denote a maximal collection of maximal subgroups of $G$ containing $H$ which are not conjugate.
Let $K_{ij}, 1 \le i \le r, 1 \le j \le n_i$ denote
the set of all maximal subgroups of $G$ containing $H$ where $K_{ij}$ is conjugate to $K_i$.

It is easy to see that $G=K_iK_j$ for $i \ne j$ (cf \cite{AsbG}).    Thus, $\Hom(\ell_0(K_i), \ell_0(K_j)=0$
if $i \ne j$.   If $K_i$ is normal in $G$, set $V_i=0$.  If not, let $V_i$ be a nontrivial irreducible submodule of $\ell_0(K_i)$
such that $\ell_0(K_i) \oplus V_i$ embeds in $\ell_0(H)$ (as in the proof of Lemma \ref{prim}).
Thus $X:= \oplus_i (W_i \oplus V_i)$ embeds
in $\ell_0(H)$ and as above, we can choose $v_{ij}$ in $W_i \oplus V_i$ linearly
independent with $K_{ij}$ the stabilizer of $v_{ij}$.
\qed

Of course, this gives:

\begin{theorem} \label{relativeth} Let $G$ be a finite solvable group.  Let $H$ be a subgroup of $G$.
Then the number of maximal subgroups of $G$ which contain $H$ is less than
$|H/G \backslash H|$.
\end{theorem}

\section{Kac algebras of Izumi-Kosaki type for solvable
groups}\label{iksection} In this section we will prove Conjecture
\ref{wall} for Kac algebras of  Izumi-Kosaki type for solvable
groups. These Kac algebras are introduced in \cite{IK} and in more
details in \cite{IK2} by considering compositions of group type
subfactors. Let us first recall some definitions from \cite{IK} to
set up our notations. The reader is refereed to \cite{IK2} for more
details. Let $G=N\rtimes H$ be semidirect product of two finite
groups $N,H.$ For $n\in N, h\in H,$ we define $n^h:=h^{-1}nh.$
Denote by $L(N)$ the set of complex valued functions on $N.$ For
$f\in L(N), f^h(n):=f(h^{-1}nh), h\in H.$
\begin{definition}\label{cocycle1}
Denote by $\eta_h(n_1,n_2), \xi_n(h_1,h_2)$ $U(1)$ valued cocycles
as defined in \S2 of \cite{IK} which verify the following cocycle
conditions:
$$
\eta_h(n_1,n_2)\eta_h(n_1n_2,n_3)=\eta_h(n_1,n_2n_3)\eta_h(n_2,n_3),
\xi_n(h_1h_2,h_3)\xi_n(h_1,h_2)=\xi_n(h_1,h_2h_3)\xi_n(h_2,h_3);
$$
Moreover, these cocycles verify the following Pentagon equation:
$$
\frac{\eta_{h_1}(n_1,n_2)\eta_{h_2}(n_1^{h_1},n_2^{h_2})}{\eta_{h_1}(n_1,n_2)}
=\frac{\xi_{n_1n_2}(h_1,h_2)}{\xi_{n_1}(h_1,h_2)\xi_{n_2}(h_1,h_2)}
$$
and normalizations:
$$
\eta_h(e,n_2)=\eta_h(n_1,e)=\xi_n(e,h_2)=\xi_n(h_1,e)=\eta_e(n_1,n_2)=1.
$$

\end{definition}
For subfactor motivations for introducing these cocycles, we refer
the reader to \S2 of \cite{IK}.
\begin{definition}\label{ikalg}
Kac algebras of Izumi-Kosaki  type are defined as Hopf algebras
$\A=L(N)\rtimes_\xi H$ whose Hopf algebra structures are given in
\cite{IK} as follows:\par (1) Algebra products:
$$
(f_1(n),h_1)(f_2(n),h_2)=(f_1(n)f_2^{h_1}(n)\xi_n(h_1,h_2),h_1h_2)
$$
where $f_2^{h_1}(n):=f_2(h_1^{-1}nh_1);$
\par
(2) Coproducts:
$$ \Delta
(n,h)=\sum_{n_2}\eta_h(nn_2^{-1},n_2)(nn_2^{-1},h)\otimes(n_2,h)
$$
\par
(3) $*$ structure:
$$
(f,h)^*= (\overline{f\xi(h,h^{-1})}^{h^{-1}}, h^{-1})
$$
\end{definition}
The following two operators on $L(N)$ will play an important role:
\begin{definition}\label{leftright}
$$
(L_{n,\eta_h}f)(m):=f(nm)\eta_h(n,m),
(R_{n,\eta_h}f)(m):=f(mn)\eta_h(n,m), \forall n,m\in N,h\in H.
$$
\end{definition}
The following lemma summarize the properties of these operators
which follow from definitions:
\begin{lemma}\label{lrprop}
$$
L_{n_1,\eta_h}L_{n_2,\eta_h}=L_{n_2n_1,\eta_h}\eta_h(n_2,n_1),
R_{n_1,\eta_h}R_{n_2,\eta_h}=R_{n_1n_2,\eta_h}\eta_h(n_1,n_2),
L_{n_1,\eta_h}R_{n_2,\eta_h}=R_{n_2,\eta_h}L_{n_1,\eta_h}.
$$
\end{lemma}
The subfactor associated with $\A$ is of the form $\L^\A\subset \L$
where $\L^\A$ is the fixed point subfactor of a factor $\L$ under
the action of $\A$ as defined in \S4 of \cite{ILP}. By \cite{ILP},
any intermediate subfactor of $\L^\A\subset \L$ is of the form
$\L^\B\subset L,$ where $\B$ is a right coideal of $\A,$ i.e., an
$*$ subalgebra of $\A$ which is verifies that $\Delta (\B)\subset
\B\otimes \A.$

The following theorem gives a characterization of  coideals of $\A$:
\begin{theorem}\label{coideal}
Let $\B$ be a right (resp. left) coideal of $\A$. Then there are
subgroups $H_1\leq H, N_1\leq N$ and $U(1)$ valued function
$\lambda:N_1\times H_1\rightarrow U(1)$ such that:\par (1) $\forall
h\in H_1, hN_1h^{-1}=N_1;$\par
\begin{equation}\label{lambdaeq}
\lambda(n_1,h)\lambda(n_2,h)=\lambda(n_1n_2,h)\eta_h(n_1,n_2),
\lambda(n,h_1)\lambda(h_1^{-1}nh_1,h)\xi_n(h_1,h_2)=\lambda(n,h_1h_2)
\end{equation}
and $\B=\oplus_{h\in H_1}(C(h),h)$ where each $C(h)\subset L(N)$
consists of functions $f\in L(N)$ such that
$L_{n,\eta_h}f=\lambda(n,h)f$ (resp. $R_{n,\eta_h}f=\lambda(n,h)f$)
$\forall n\in N_1, h\in H_1$.\par Conversely, any triple
$(N_1,H_1,\lambda)$ which verify the above conditions uniquely
determine a coideal of $\A.$
\end{theorem}
\prf We will prove the theorem for the case when $\B$ is a right
coideal of $\A$. The remaining case is similar. We write elements of
$\B$ as $\sum_h (f_h,h)$ where $f_h\in L(N).$ We have
$$
\Delta \sum_h(f_h,h)= \sum_{n_2,h}
(R_{n_2,\eta_h}f_h,h)\otimes(n_2,h)
$$
Since $\B$ is a right coideal, it follows that for each fixed
$(n_2,h)$, $(R_{n_2,\eta_h}f_h,h)\in \B.$ So we have
$\B=\oplus_h(C(h),h)$ with $C(h)$ a subspace of $L(N)$ which is
mapped by $R_{n,h}$ to itself. Since $\B$ is also an algebra, we
have
\begin{equation}\label{cs}
C(h_1)C(h_2)^{h_1}\xi(h_1,h_2)\subset C(h_1h_2)\\
\end{equation}
In particular $C(e)$ is a subalgebra of $L(N)$ which affords a right
representation of $N.$ It follows that there is a subgroup $N_1\leq
N$ such that $C(e)$ is the space of $N_1$-left invariant functions
on $N.$ Let $N=\bigcup_i N_1b_i, 1\leq i\leq k$ with $k=|N|/|N_1|$
be the left coset decompositions of $N.$ Then $\delta_{Nb_i}$ is a
basis of $C(e).$

\par
 Since $\B$ is an $*$ algebra, it follows that if
$(f_h,h)\in \B,$ then
$$(f_h,h)^*\in
(C(h^{-1}),h^{-1}),(f_h,h)(g_h,h)^*= (f_h\bar{g_h},e)\in (C(e),e).$$
Let $H_1:=\{h\in H| C(h)\neq 0\}.$ It follows easily from above that
$H_1$ is a subgroup of $H.$ By equation (\ref{cs}) $C(e)C(h)\subset
C(h).$ so it follows that $C(h)=\oplus_i C(h)\delta_{N_1b_i}.$
Assume that $h\in H_1$ so $C(h)\neq 0.$ Since $R_{b_i^{1},\eta_h}$
maps $C(h)$ to itself, we can assume that $C(h)\delta_{N_1}\neq 0.$
Let $f_i\neq 0\in C(h)\delta_{N},i=1,2$ then $f_1\bar f_2\in
C(e)\delta_{N_1}=\C\delta_{N_1}.$ We conclude that
$C(h)\delta_{N_1}$ is one dimensional, and by using operator
$R_{b_i^{-1},\eta_h},$ we conclude that $C(h)\delta_{N_1b_i}$ is one
dimensional. Let $f_h\neq 0\in C(h)\delta_{N_1},$ then since
$R_{n_1,\eta_h}$ maps $C(h)\delta_{N_1}$ to itself, there is a
function $\lambda:N_1\times H_1\rightarrow U(1)$ such that
$$
R_{n_1,\eta_h}f_h=\lambda(n_1,h)f_h.
$$
We can assume that $f_h(e)=1.$ Then we have
$$
(R_{n_1,\eta_h}f_h)(e)=\lambda(n_1,h)=f_h(n_1).
$$
Equation (\ref{lambdaeq}) follows from Lemma \ref{lrprop} and
equation (\ref{cs}). Let us show that $C(h),h\in H_1$ is the
subspace of $L(N)$ which verifies
$$
L_{n_1,\eta_h}f=\lambda(n_1,h)f.
$$
We note that by definition \ref{leftright} $f_h(n_1)=\lambda(n_1,h)$
verifies
$$
(L_{n_1,\eta_h}f_h)(m)=\lambda(n_1m,h)\eta_h(n_1,m)=\lambda(n,h)\lambda(m,h)
=\lambda(n,h)f_h(m)
$$
where in the last equation we have used equation (\ref{lambdaeq}).
Since $C(h)$ is the linear span of $R_{n,\eta_h}f_h$, by lemma
\ref{lrprop} we have proved that for any $f\in C(h),
L_{n_1,\eta_h}f=\lambda(n_1,h)f.$ On the other hand by counting
dimensions we conclude that $C(h),h\in H_1$ is the subspace of
$L(N)$ which verifies
$$
L_{n_1,\eta_h}f=\lambda(n_1,h)f.
$$
Let us show that $\forall h\in H_1, hN_1h^{-1}=N_1.$ By equation
(\ref{cs}) we have $f_h(n_1)\delta_{N_1}^h\subset C(h)\delta_{N_1}$
which is one dimensional, and it follows that
$\delta_{N_1}^h=\delta_{N_1}$ for all $h\in H_1,$ i.e.,
$hN_1h^{-1}=N_1.$\par

Conversely for  any triple $(N_1,H_1,\lambda)$ which verify the
conditions in theorem \ref{coideal}, we can simply define
$\B:=\oplus_{h_1\in H_1}(C(h),h)$ where $C(h)\subset L(N)$ consists
of functions $f\in L(N)$ such that $L_{n,\eta_h}f=\lambda(n,h)f.$ We
need to check that $\B$ is a right coideal. By inspection it is
enough to check equation (\ref{cs}). By definition we need to check
that if $f_{h_i}\in C(h_i), i=1,2,$ then
$g(n):=f_{h_1}(n)f_{h_2}^{h_1}(n)\xi_n(h_1,h_2)\in C(h_1h_2).$ So we
need to show that $(L_{n,\eta_{h_1h_2}}g)(m)=\lambda(n,h_1h_2)g(m),
\forall n\in N_1, m\in N.$ Since $f_{h_i}\in C(h_i), i=1,2,$ we have
$$
f_{h_i}(nm)\eta_{h_i}(n,m)=\lambda(n,h_1)f_{h_i}(m),i=1,2.
$$
By using above equation and equation (\ref{lambdaeq}) it follows
that $(L_{n,\eta_{h_1h_2}}g)(m)=\lambda(n,h_1h_2)g(m), \forall n\in
N_1, m\in N$ iff the following holds:
$$
\frac{\eta_{h_1}(n,m)\eta_{h_2}(n^{h_1},m^{h_2})}{\eta_{h_1h_2}(n,m)}
=\frac{\xi_{nm}(h_1,h_2)}{\xi_{n}(h_1,h_2)\xi_{m}(h_1,h_2)}
$$

which is the pentagon equation in definition \ref{cocycle1}.

\qed

For a coideal $\B$ with $(N_1,H_1,\lambda)$ as in Th. \ref{coideal},
we shall refer to $(N_1,H_1,\lambda)$ as {\it the triple} associated
with $\B.$ We note that by Th. \ref{coideal}, such triple uniquely
determine $\B.$ Moreover, suppose that the triples associated with
$\B_i$ are given by $(N_i,H_i,\lambda_i), i=1,2.$ Then $\B_1\subset
\B_2$ iff $N_1\supset N_2, H_1\subset H_2$ and $\lambda_1,\lambda_2$
agree on $N_2\times H_1.$
\begin{lemma}\label{count1}
Let $\B$ be a right coideal of $\A$ as in Theorem \ref{coideal} with
triple $(N_1,H_1, \lambda).$ The the number of right coideals of
$\A$ with the same $(N_1,H_1)$ are given as follows: Let $\hat{N_1}$
be the set of homomorphisms from $N_1$ to $U(1)$ and form a group
$\hat{N_1}\rtimes H_1$. Then the  right coideals of $\A$ with the
same $(N_1,H_1)$ are in one to one correspondence with the set of
cocycles from $H_1$ to $\hat{N_1},$ i.e., maps $\mu: H_1\rightarrow
\hat{N_1}$ such that
$$
\mu(h_1)\mu(h_2)^{h_1}=\mu(h_1h_2).
$$
\end{lemma}
\prf Let $\B_1$ be a right coideal of $\A$ with triple
$(N_1,H_1,\lambda_1)$ as in Theorem \ref{coideal}. Let
$\mu:=\lambda_1/\lambda.$ By equation (\ref{lambdaeq}) we conclude
that $\mu$ is a cocycle from $H_1$ to $\hat{N_1}.$ Conversely, if
$\mu$ is a cocycle from $H_1$ to $\hat{N_1},$ then $\B_1$ associated
with the triple $(N_1,H_1,\lambda \mu)$ is a  right coideal of $\A$
by Theorem \ref{coideal}. \qed
\begin{lemma}\label{count2}
Let $H_1\subset H, N_1\neq \{e\}\subset N$ such that
$hN_1h^{-1}=N_1, \forall h\in H_1.$ Let $\hat{N_1}$ be the set of
homomorphisms from $N_1$ to $U(1)$ and form a group
$\hat{N_1}\rtimes H_1$. Assume that $H_1$ acts irreducibly on
$\hat{N_1}$ and $H_1$ is solvable. Then the number of cocycles from
$H_1$ to $\hat{N_1},$ i.e., maps $\mu: H_1\rightarrow \hat{N_1}$
such that
$$
\mu(h_1)\mu(h_2)^{h_1}=\mu(h_1h_2)
$$
is less or equal to  $(|N_1|-1)|H_1|.$
\end{lemma}
\prf If $H_1$ acts trivially on $\hat{N_1},$ since $H_1$ acts
irreducibly on $\hat{N_1},$ it follows that $\hat{N_1}$ is an
abelian group of prime order, and the number of cocycles from $H_1$
to $\hat{N_1}$ is bounded by $|H_1|.$ If  $H_1$ acts nontrivially on
$\hat{N_1},$ then $H_1$ is a maximal subgroup of $\hat{N_1}\rtimes
H_1$ with trivial core, and since $H_1$ is solvable, all cocycles
from $H_1$ to $\hat{N_1}$ are coboundaries by Th.16.1 of
\cite{solvable}, and is bounded by $|\hat{N_1}|.$ Since $|H_1|>1$ in
this case the lemma is proved. \qed
\begin{theorem}\label{depth2th}
Let $\A=L(N)\rtimes_\xi H$ be Kac algebras of Izumi-Kosaki type as
in definition \ref{ikalg}. Assume that $N,H$ are solvable groups.
Then the number of maximal (resp. minimal) right coideals is less
than the dimension of $\A.$
\end{theorem}
\prf Assume that $\B$ is a right coideal of $\A$ and let
$(N_1,H_1,\lambda)$ be the triple as in Th. \ref{coideal}. Let us
first prove the minimal case. Since $\B\supset (L(N/N_1),e),$ if
$N_1\neq N,$ we must have $\B= (L(N/N_1),e),$ and $N_1$ must be
maximal in $N.$ The number of such $N_1$ is less than $|N|$ by Th.
\ref{relativeth}.\par If $N_1=N,$ then $H_1\neq e.$ Let $\Z_p\leq
H_1$ be any minimal subgroups of $H_1,$ then the triple
$(N,\Z_p,\lambda)$ will give rise to a right coideal of $\A$ by Th.
\ref{coideal} which is contained in $\B.$ It follows that
$H_1=\Z_p.$ By Lemma \ref{count1}, the number of such triple is
bounded by $|\hat{N}|\leq |N|.$ So minimal right coideal of $\A$ is
bounded by the sum of number of maximal subgroups of N and the
product of the number of minimal subgroups of $H$ by $|N|,$ and it
follows that the number of minimal right coideals is less than the
dimension of $\A.$\par Now assume that $\B$ is maximal. If $N_1$ is
trivial, then $H_1$ is maximal in $H,$ and by Th. \ref{relativeth}
the number of  maximal $H_1$ is less than $|H|.$\par If $N_1$ is
nontrivial, then $\bigoplus_{h\in H_1}(L(N),h)\supset \B,$ and it
follows that $H_1=H.$ We claim that $N_1$ is generated by $\Ad_H
(x)$ for any nontrivial $x\in N_1.$ In fact let $N_1'\subset N_1$ be
a subgroup generated by  $\Ad_H (x)$ for a nontrivial $x\in N_1.$
Then the right coideal determined by the triple $(N_1',H,\lambda)$
contains $\B$, and by maximality of $\B$ we have $N_1'=N_1.$ It
follows that $H$ acts irreducibly on $N_1/[N_1,N_1],$ and therefore
acts irreducibly on its dual $\hat{N_1}.$ By lemma \ref{count2} such
$\B$ with fixed $(N_1,H)$ is bounded by $(|N_1|-1)|H|.$ Note that
different $N_1$'s intersect only at identity. It follows that the
number of maximal $\B$'s is bounded by
$$
(|H|-1)+|H|(|N|-1)=|H||N|-1.
$$
\qed

We consider the intermediate subfactors of $\L^\B\subset L$
corresponding to $\B$ as in Th. \ref{coideal}.
\begin{lemma}\label{jb}
The dimension of second higher relative commutant associated with
the subfactor $\L^{\B}\subset \L$ is given by
$$
\sum_{h\in H_1} \dim(C_R(h)\cap C_L(h))
$$
where
$$
C_R(h)\cap C_L(h):=\{f\in L(N),
R_{n_1,\eta_h}f=\lambda(n_1,h)f=L_{n_1,\eta_h}f\}.
$$
\end{lemma}
\prf This follows from \S3 of \cite{NV} and  Th. \ref{coideal}.\qed
\begin{lemma}\label{jb1}
Let $\B$ be as in Th. \ref{coideal} with triple $(N_1,H_1,\lambda).$
Suppose that $\lambda$ can be extended to $N_i\supset N_1$ such that
the triple $(N_i,H_1,\lambda),i=1,2,...,n$ gives a right coideal of
$\A$ via Th. \ref{coideal}, and $N_i\cap N_j=N_1, \forall i\neq j.$
Let $k_i$ be the number of double cosets of $N_1$ in $N_i.$ Then
$$\dim(C_R(h)\cap C_L(h))\geq 1+(k_2-1)+...+(k_n-1).
$$
\end{lemma}
 \prf On each double coset $N_1bN_1$ of $N_i,$ we can define a
function such that its value on the double coset is simply the value
of $\lambda (.,h)$ and zero elsewhere. It is easy to check that
these functions belong to $C_R(h)\cap C_L(h),$ and they are linearly
independent since they have different support, and the lemma
follows. \qed

The following two lemmas are straightforward consequences of
definitions:
\begin{lemma}\label{db1}
Let $N_2\supset N_1.$ Then the number of homomorphisms from $N_2$ to
$U(1)$ which  takes value $1$ on $N_1$ is bounded by the number of
double cosets of $N_1$ in $N_2.$
\end{lemma}
\begin{lemma}\label{miniext}
Let $N_2\subset N$ be a minimal extension of $N_1$ which is
$\Ad_{H_1}$ invariant. Then the natural action of $H_1$ on
$N_2/N_1[N_2,N_2]$ is irreducible. The dual of $N_2/N_1[N_2,N_2]$ is
abelian group of homomorphisms from $N_2$ to $U(1)$ which  takes
value $1$ on $N_1.$
\end{lemma}
\begin{theorem}\label{relative2}
Let $\B\subset \A$ be a right coideal of $\A$ as in Th.
\ref{coideal}. Then both minimal version and maximal version of
Conjecture \ref{wall} is true for subfactors $\L^\B\subset \L$ when
$H,N$ are solvable.
\end{theorem}
\prf Let $\B_i\subset \B$ and $(N_i,H_i,\lambda_i)$ be the
associated triple of $\B_i$ as in Th. \ref{coideal}. We have
$N_i\supset N_1, H_i\subset H_1,$ and $\lambda_i$ agrees with
$\lambda$ on $N_1\times H_i.$\par Let us first prove the maximal
case. Assume that $\B_i\subset \B,i=1,2,3,..,m$ is the list of
maximal right coideals of $\A$ that is contained in $\B$. Let
$(N_j,H_j), j=1,2,3,...,n$ be the list of different pairs of groups
that are associated with $\B_i's$ as in Th. \ref{coideal}.

If $N_j=N_1,$ then $H_j\subset H_1$ is maximal in $H_1.$ The number
of such maximal $H_j$ is less than $|H_1|$ by Th. \ref{relativeth}.
If $N_j\neq N_1,$ since $\Ad_{H_j}(N_1)=N_1,$ it follows that
$H_j=H_1,$ and $N_j$ is a minimal extension of $N_1$ which is
invariant under $\Ad_{H_1}.$ Let $k_j$ be the number of double
cosets of $N_1$ in $N_j.$ Note that these $N_j$'s only intersect at
$N_1. $ By Lemma \ref{miniext},Lemma \ref{db1} and Lemma
\ref{count2} the number of such $\B_i$ with fixed $N_j$ is bounded
by $(k_j-1)|H_1|.$ So $m\leq |H_1|-1 + [(k_2-1)+...+(k_n-1)]|H_1|$
and the theorem follows from Lemma \ref{jb} and Lemma \ref{jb1}.\par
Now assume that $\B_i\subset \B,i=1,2,3,..,m$ is the list of minimal
right coideals of $\A $ that is contained in $\B$. Let $(N_j,H_j),
j=1,2,3,...,p$ be the list of different pairs of groups that are
associated with $\B_i's$ as in Th. \ref{coideal}. By considering
$(L(N/N_j),e)$ it follows that if $H_j$ is trivial, then $N_j\subset
N$ is a maximal subgroup, and the number of such maximal subgroups
is bounded by the double cosets of $N_1$ in $N$ by Th.
\ref{relativeth}.  If $H_j$ is nontrivial, then $N_j=N,$ and it
follows that $H_j$ has to be a minimal nontrivial subgroup of $H_1.$
The number of such subgroups of $H_1$ is bounded by $|H_1|-1$ if
$H_1$ is not an abelian group of prime order, and $1$ if $H_1$ is an
abelian group of prime order. For each fixed $(N,H_j),$ the possible
$\lambda_j$'s which agrees with $\lambda$ on $N_1\times H_j$ is
clearly bounded by the number of homomorphisms from $N$ to $U(1)$
which vanishes on $N_1$, and by Lemma \ref{db1} this number is
bounded by the number of double cosets of $N_1$ in $N$ which is
denoted by $p_1$.  It follows that $p\leq p_1-1+ (|H_1|-1)p_1,$ and
by Lemma \ref{jb} we are done.\qed
\begin{remark}
If we set $H$ to be a trivial group in Th. \ref{relative2}, then we
recover Th. \ref{relativeth}.
\end{remark}

\begin{remark}
For subfactor $\L^\A\subset \L^\B,$ we can map its intermediate
subfactors to certain right coideals of $\hat{\A}$ which is the dual
of $\A$ (cf. \S2 of \cite{IK}). Essentially similar argument as in
the proof of Th. \ref{relative2} shows that subfactor
 $L^\A\subset L^\B$ verifies the maximal and minimal version of Conjecture \ref{wall} when $H,N$ is
solvable.
\end{remark}

\section{Tensor Product Conjecture}\label{tensorsection}
\begin{lemma}\label{tensorl}
Let $G = X \times Y$ be finite groups with $|X|=x$and $|Y|=y$. Then
the number  of maximal subgroups of G which contain neither $X$ nor
$Y$ is at most $(x-1)(y-1)$ (with equality if and only if X and Y
are elementary abelian 2-groups).
\end{lemma}

\prf  Let $M $ be a maximal subgroup of $G$ containing neither $X$
nor $Y$. Let  $f~:~G\rightarrow G/K$ be the natural homomorphism where
$K$ is the core of $M $ in $G$. Then $f(X)$ and $f(Y)$ are normal
nontrivial subgroups which commute in the primitive group $G/K$ and
moreover, they generate $G/K$ together. By the
Aschbacher-O'Nan-Scott Theorem (although this can be proved easily
in this case) this implies that either $f(X)=f(Y)$ has prime order
$p$ for some prime $p$ or $G/K = f(X) \times f(Y) = S \times S$ with
$S$ a nonabelian simple group.

Thus,  passing to the quotient  by the intersection of all the cores of such
maximal subgroups, we may assume that $X$ and $Y$ are direct products
of simple groups.  Write $X=  \prod_p X_p \times \prod_S X_S$ where $X_p$
is the maximal elementary abelian $p$-quotient of $X$ and $X_S$
is the maximal quotient of $X$ that is a direct product of nonabelian simple groups
each isomorphic to $S$.   Write $Y$ in a similar manner.  The previous paragraph
shows that we may reduce to the case that either $X$ and $Y$ are each elementary
abelian $p$-groups or are both direct products of a fixed nonabelian simple group $S$.

In the first case, it is trivial to see that the total number of maximal subgroups
is $(xy-1)/(p-1)$ while the number of maximal subgroups containing either $X$ or $Y$
is  $(x-1)/(p-1) + (y-1)/(p-1)$.  Thus, the total number of maximal subgroups containing
neither $X$ nor $Y$ is $(x-1)(y-1)/(p-1)$.

In the second case,  write $X=S^a$ and $Y=S^b$.  If $M$ is a maximal subgroup
not containing $X$ or $Y$, then $M \cap X$ is normal in $X$ with
$X/(M \cap X) \cong S$.   Thus $M$ is a direct factor of $X$
isomorphic to $S^{a-1}$.  There are $a$ such factors.    Thus, the number
of maximal subgroups of $X \times Y$ not containing $X$ or $Y$ is $abc$
where $c$ is the number of maximal subgroups of $S \times S$ not containing
either factor.  This is precisely $|\aut(S)|$  (since any such maximal subgroup
is $\{s, \sigma(s) | s \in S\}$ where $\sigma \in \aut(S)$).   Thus, in this case
the number of maximal subgroups not containing either factor is $ab|\aut(S)|$.
To complete the proof, we only need to know that $|\aut(S)| <( |S|-1)^2$.

This is well known (and in fact much better inequalities can be shown).   All such
existing proofs depend upon the classification of finite simple groups.  We note
that the inequality we need follows from the fact that every finite nonabelian
simple group can be generated by two elements  (note that if $s \in S$, then
$|\{\sigma(s)| \sigma \in \aut(S)\} < |S| -3 $ (for there are at least $4$ different
orders of elements and so at the very least $4$ different orbits on $S$  -- if
$s,t$ are generators, then any automorphism is determined by its images
on $s,t$ whence the inequality).
\qed

The following corollary follows immediately:

\begin{corollary}\label{tensorw}
Let $G = X \times  Y$ be finite groups such that both $X$ and $Y$
verify Wall's conjecture, then $G$ also verifies Wall's conjecture.
\end{corollary}
Based on Lemma \ref{tensorl}, we propose the following tensor
product conjecture:
\begin{conjecture}\label{tensorc}
Let $N_i\subset M_i, i=1,2$ be two irreducible subfactors with
finite index. Then the number of minimal intermediate subfactors in
$N_1\bigotimes N_2\subset M_1\bigotimes M_2$ which is not of the
form $N_1\bigotimes P, P\bigotimes N_2$ is less or equal to
$(n_1-1)(n_2-1)$ where $n_i$ is the dimension of second higher
relative commutant of $N_i\subset M_i, i=1,2.$
\end{conjecture}
This conjecture is nontrivial even for subfactors coming from
groups, where we have seen in the proof of Lemma \ref{tensorl} that
we have used classification of finite simple groups to bound the
number of automorphisms of a simple group.

\end{document}